\documentclass[11pt]{article}
\usepackage[all]{xy}
\usepackage{graphicx}
\usepackage[margin=2.95cm]{geometry}
\usepackage{amssymb}
\usepackage{latexsym} 
\newtheorem{thm}{Theorem}[section]

\newtheorem{prop}[thm]{Proposition} 
 \newtheorem{dfn}[thm]{Definition}

 \newtheorem{rmk}[thm]{Remark}
\newtheorem{ex}[thm]{Example} 

\newcommand{\pf}{\noindent{\bf Proof.}\ }

\newcommand{\Hom}{{\rm Hom}}

\newcommand{\im}{{\rm Im \,}}
\newcommand{\dom}{{\rm Dom \,}}

\newcommand{\xbar}{{\overline{X}}}
\newcommand{\ybar}{{\overline{Y}}}

\newcommand{\boldc}{{\mathbf C}}

\newcommand{\boldlcf}{{\mathbf f}}
\newcommand{\boldlcg}{{\mathbf g}}

\newcommand{\qed}{\begin{flushright} $\Box$\ \ \ \ \ \end{flushright}}

\newcommand{\from}{\leftarrow}
\newcommand{\mapsfrom}{\reflectbox{$\mapsto$}}
\newcommand{\rel}{\mathbf{REL}}

\newcommand{\mrel}{\mathbf{MREL}}
\newcommand{\srel}{\mathbf{SREL}}

\newcommand{\unitob}{\mathbf{1}}

\newcommand{\arrows}{\,\lower1pt\hbox{$\longrightarrow$}\hskip-.24in\raise2pt
             \hbox{$\longrightarrow$}\,}

\newcommand{\tr}{\mathrm{t}}

\newcommand{\xfy}{X \stackrel{f}{\longleftarrow}Y}

\newcommand{\ygz}{Y \stackrel{g}{\longleftarrow}Z}

\newcommand{\xfgz}{X \stackrel{f\circ g}{\longleftarrow}Z}

\begin{document}

\message{ !name(hyperlagrangians.tex) !offset(-3) }

\title{{\bf A note on the Wehrheim-Woodward category}}
\author
{Alan
Weinstein
\thanks{Research partially supported by NSF Grant
DMS-0707137.
\newline \mbox{~~~~}MSC2010 Subject Classification Number: 
53D12 (Primary), 18B10 (Secondary).
\newline \mbox{~~~~}Keywords: symplectic manifold, canonical relation,
lagrangian submanifold, categories of relations}
\\
Department of Mathematics\\ University of California\\Berkeley, CA
94720 USA\\ {\small(alanw@math.berkeley.edu)}}
\date{}
\maketitle

\begin{center}
{\em Dedicated to Tudor Ratiu for his 60th birthday}
\end{center}

\bigskip

\begin{abstract} 
Wehrheim and Woodward have shown how to embed all the canonical
relations between symplectic manifolds into a
category in which the composition is the usual one when transversality and
embedding assumptions are satisfied.   A morphism in their category is
an equivalence class of composable sequences of canonical relations,
with composition given by
concatenation.  In this note, we show that every such morphism is
represented by a sequence consisting of just two relations, one of them a
reduction and the other a coreduction.
\end{abstract}


\section{Introduction}
\label{sec-intro}
The problem of quantization, i.e., the transition from
classical to quantum physics, may be formulated mathematically as the
search for a functor from a ``classical'' category whose
objects are symplectic manifolds to a ``quantum'' category whose objects are Hilbert
spaces (or more general objects, such as spaces of distributions, Fukaya categories or
categories of D-modules).  

On the classical side, it is useful to include among the morphisms 
$X\from Y$ not
only symplectomorphisms, which should produce unitary operators upon
quantization, but more general canonical relations, i.e. lagrangian
submanifolds of $X\times\ybar$ (where $\ybar$ is $Y$ with its
symplectic structure multiplied by $-1$).  An immediate difficulty is that the
composition of canonical relations can produce relations which are not even
smooth submanifolds. 

A solution to the composition problem on the classical side has been
given by Wehrheim and Woodward \cite{we-wo:functoriality}, 
who introduce a category which is
generated by the canonical relations, but where composition is merely
``symbolic'' unless the pair being composed fits together in the best
possible sense.  In other words, the morphisms in this category are
equivalence classes of sequences of canonical relations which are
composable as set-theoretic relations.  The equivalence relation
allows one to shorten a sequence when two adjacent entries compose
nicely.

The purpose of this note is to show that 
any morphism in the WW category may be expressed as a product
of just {\em two} canonical relations.
Furthermore, this
factorization $g\circ h$ is analogous to the factorization of
a ordinary map $X\from Y$ through $X\times Y$ as the product of
the projection $X\from X\times Y$ and the graph map $X\times
Y\from Y$ in that the canonical relation $g$ is the
symplectic analogue of a submersion (it is a kind of symplectic
reduction),  and $h$ 
is analogous to an embedding.   In addition, the natural transpose
operation on relations, which exchanges source and target, extends to the
WW category, and the subcategories to which $g$ and $h$ belong are
transposes of one another; in this sense, relations have an even nicer
structure than maps.

The decomposition $g\circ h$ is not unique; in the iterative 
construction which we will
describe, the source of $g$ and target of
$h$ is a product of spaces whose number
 grows exponentially with $n$.  Another, 
due to Katrin
 Wehrheim, consists of just $n$ or $n+1$ factors, depending on the
 parity of $n$.  

In a longer paper in preparation, we will show that the
Wehrheim-Wood\-ward construction leads to a rigid monoidal category,
in which each morphism $\xfy$ may be represented by a ``hypergraph'' which
is a lagrangian submanifold of a symplectic manifold $Q$ of which
$X\times \ybar$ is a symplectic reduction.  We will also place the 
Wehrheim-Woodward construction in an appropriate general setting 
which is robust enough to apply to 
categories of operators in which quantization
functors may take their values. 

\bigskip
\noindent
{\bf Acknowledgements.}  I would like to thank the Institut
Math\'ematique de Jussieu for many years of providing a stimulating
environment for research during my annual visits.  
For helpful comments, I would like to thank Sylvain Cappell, Thomas
Kragh, Sikimeti Ma'u, Katrin
Wehrheim, and Chris Woodward, as well as the referees.

\section{Relations and their composition}
\label{sec-relations}
This section is a review of mostly well-known ideas concerning the
category of sets and relations, with some new terminology and notation.

We denote by $\rel$ the category whose objects are sets
and for which the morphism space $\rel(X,Y)$ is the
set of all subsets of $X\times Y$.   We adopt the convention that
$f \in \rel(X,Y)$ is a morphism \emph{to} $X$ \emph{from} $Y$.
Thus, $X$  is the {\bf target} of $f$ and $Y$
the {\bf source}, and we write $\xfy$ for this morphism.
The {\bf composition} $\xfgz$ of $f\in \rel(X,Y)$ and $g\in \rel(Y,Z)$
is $$\{(x,z)| (x,y)\in f ~{\rm and}~
(y,z) \in g ~{\rm for~some}~y\in Y\}.$$   
We denote the identity relation on $X$ by
$X\stackrel{1_X}{\longleftarrow} X$,
but also by $\Delta_X$ when we want to think of it explicitly as a
subset of $X\times X$.  

The natural exchange maps $X\times Y \leftarrow Y\times X$ define
an involutive
contravariant {\bf
  transposition} functor $f\mapsfrom f^\tr$
 from $\rel$ to itself.  It is the identity on objects.

For any subset $T$ of $Y$, the {\bf image} $f(T)$
is the subset  $\{x\in X|(x,y)\in f~ {\rm for~some}~ y\in T\}$ of $X$.
For a single element $y\in Y$, we write $f(y)$ for the subset $f(\{y\})$.
 The image $f(Y)$ is called the
{\bf range} of $f$, and the range  $f^\tr(X) \subseteq Y$ of the
transpose is the {\bf  domain}.  We denote these by $\im f$ and $\dom
f$ respectively.

$f$ is {\bf surjective} if its range equals its target, and {\bf
  cosurjective} if its domain equals its source (i.e. if it is
``defined everywhere'').  $f$ is {\bf injective} if, for any $x\in X$,
there is at most
one $(x,y) \in f$, and {\bf coinjective} if there is at most one
$(x,y)\in f$ for any $y\in Y$ (i.e. if it is ``single valued'').  
Thus, the cosurjective and coinjective relations in $X\times Y$ are
the graphs of functions to $X$ from $Y$.  $f$ is coinjective
[cosurjective] if and only if $f^\tr$ is injective [surjective].

Each of
the four classes just defined constitutes a subcategory of $\rel$, but
two intersections of these subcategories
 will be particularly important.  The surjective and coinjective relations 
 (i.e. partially defined but 
single valued surjections) will
be called {\bf reductions}, and those which are injective and
cosurjective (i.e. which take all the  points of the source to disjoint
nonempty subsets of the target) will be called {\bf coreductions}.
(In the linear context, these special relations are called reductions
and contrareductions by Benenti and Tulczyjew \cite{be-tu:relazioni}, who
characterize them by several categorical properties.)
We will sometimes indicate that a morphism is of one of these types by
decorating the arrow which represents it: 
$X \twoheadleftarrow Y$ for a reduction and $X\leftarrowtail Y$ for a
coreduction.   The arrows decorated at both ends are the invertible relations, i.e., the bijective functions.

The composition of relations 
$\xfy$ and $\ygz$ involves  two steps.  
The first is to form the fibre product $f \times_Yg $, i.e. the intersection 
 in $X\times Y\times
Y\times Z$ of $f \times g$ with $X \times
\Delta_Y \times Z$.  The second is to  project $f\times_Y g$
into $X \times Z$.  We will call the pair
$(f,g)$ {\bf monic} if this projection 
is injective, i.e. 
if there is only one $y\in Y$ which accounts for 
each $(x,z)$ belonging to $f\circ g$.

A composable pair $(f,g)$ is automatically monic if $f$ is injective
or $g$ is coinjective.  Thus, all compositions within the categories
of injective relations or coinjective relations, and hence of
reductions or coreductions, are monic.  

\section{Smooth relations}
The objects on which smooth relations
operate are not just sets.  To get a category, we begin
with
all relations between manifolds before singling out the smooth ones.

\begin{dfn}
The objects of the category $\mrel$ are the smooth differentiable
manifolds, and the morphisms are the set-theoretic relations; 
i.e., $\mrel(X,Y)$ consists of all sub{\em sets} of the product
manifold $X\times Y$.   A relation $f\in \mrel(X,Y)$ is {\bf smooth}
if it is a closed submanifold of $X\times Y$.  The manifold whose only
element is the empty set, carrying its unique smooth structure, will be denoted by $\unitob$.
\end{dfn} 

The forgetful functor $\rel \from \mrel$ which forgets the
differentiable structure is full as well as faithful.

\begin{ex}
{\em
The graph of any smooth map $\xfy$ is a smooth relation.  
The smooth relations in $\mrel(Y,\unitob)$ correspond to the closed submanifolds
$M\subseteq Y$.  The composition in $\mrel$ of the relations
corresponding to $f$ and $M$ corresponds
to the subset $f(M)$, which is in general neither closed nor a
submanifold, hence not a smooth relation.
}
\end{ex}

We see already from the example above that we will need to impose
extra conditions on pairs of smooth relations to insure that their
composition is smooth.   To describe these conditions, it is useful to
look at the tangent ``operation'', which is not a functor since its
domain is not (yet) a category. 

Since the tangent bundle of any submanifold is a smooth submanifold of
the pulled back tangent bundle of the ambient manifold, 
the tangent bundle $Tf$ of any smooth relation $f\in \mrel(X,Y)$ is a 
smooth relation in
$\mrel(TX,TY)$.  We call it the {\bf differential} of $f$.

For smooth relations, we will use more restrictive definitions of
reduction and coreduction.
A smooth relation $f \in \mrel(X,Y)$ will be called a {\bf reduction} if it is
surjective and coinjective, along with $Tf$, and if the projection
$Y\from f$ is not only injective but proper. 
$f$ will be called an {\bf coreduction} if $f^t$ is a reduction.
Note that a 
reduction can be partially defined but is single valued, along with
its differential; a coreduction can be multiply
defined but, along with its differential, must be defined everywhere.   

We will also reserve the term ``monic'' for compositions 
$\xfy\stackrel{g}{\longleftarrow}Z$ where the
projection $X\times Z\from f\times_Y g$ is not only injective but
also proper.  

For compositions of smooth relations, we may impose a transversality
condition which does not seem to have a useful counterpart in $\rel$.

\begin{dfn}
\label{dfn-smoothtransversal}
When $\xfy$ and $\ygz$ are smooth,  the 
pair $(f,g)$ is {\bf transversal} if 
$f\times g$ is
transversal to $X \times
\Delta_Y \times Z$, insuring that $f\times_Y g$ is again a
manifold.   
A transversal pair is {\bf strongly transversal} if 
the projection map  $X\times Z \from f\times_Y g$ is an embedding
onto a closed submanifold (which is $f\circ g$).  This implies
that $(f,g)$ and $(Tf,Tg)$ are both monic, and a consequence of
these conditions is that $\xfgz$ is again a smooth relation.
\end{dfn}

\begin{rmk}
{\em
For smooth transversal $(f,g)$, 
the monicity condition at each $(x,y,y,z) \in f\times_Y g$
applied to $T_{(x,y)}f$ and $T_{(y,z)}g$ means that 
the projection from $f\times_Y g$
to $f\circ g$ 
is an immersion, but we need $(f,g)$ itself to be monic to insure that this
immersion is an embedding, so that $f$ is strongly transversal to $g$.    
}
\end{rmk}

Note that a composition $X\leftarrow Y\leftarrow Z$ 
is strongly transversal if either arrow is decorated at the junction
point $Y$, i.e. if 
we have either $X\leftarrowtail Y$ or $Y\twoheadleftarrow Z$ (or
both).  The  
reductions and coreductions between 
manifolds
each form
subcategories of $\mrel$.  Other subcategories are given by the graphs
of smooth maps and by the transposes of such graphs.    

\begin{ex}
{\em
If $f$ and $g$ are smooth, and their composition
$f\circ g$ happens to be a smooth relation, it may still fail to be the case that
$T(f\circ g) = Tf \circ Tg$.   

For instance, let $C$ and $D$ be submanifolds
in a manifold $X$ which intersect in a single point $x$, but whose
tangent spaces there are equal.  Let $f$ be  the reduction $C\twoheadleftarrow X$
whose domain is $C$ (transpose of the inclusion), and let $g$ be the inclusion
$X\leftarrowtail D$.  The composition $f\circ g$ is the single point
$(x,x)$ in $C \times D$, but the composition $Tf\circ Tg$ is the
graph in $TC\times TD$ of the 
isomorphism $T_x C \from T_x D$ given by their inclusion in $T_xX$.
}
\end{ex}

\section{Canonical relations}

A {\bf canonical relation} (called in \cite{we-wo:functoriality} a``lagrangian
correspondence'')
between symplectic manifolds $(X,\omega_X)$ 
and $(Y,\omega_Y)$ 
is a smooth relation $\xfy$ which
is lagrangian as a submanifold of $(X,\omega_X) \times (Y,-\omega_Y)$. 
We will often omit the symbol for the symplectic structure
and will use the notation
$\xbar$ for the {\bf dual} $(X,-\omega_X)$ when $X$ is $(X,\omega_X)$.  

If we drop the condition that $f$ be canonical, or even that it be
smooth, we obtain the category $\srel$ of  relations between
symplectic manifolds.   The canonical relations are morphisms in this
category, and the starting point for their study is the important fact
that if $\xfy$ and $\ygz$ are canonical, and if $f$ is strongly
transversal to $g$, so that $f \circ g$ is smooth, then
$\xfgz$ is canonical as well.   The graph $\gamma_f$ of $\xfy$ 
will be considered as an element of $\srel(X\times
\ybar,\unitob)$.   

It is a well-known and useful fact that any smooth map $\xfy$
may be factored as the composition $g\circ h$ of a surjective submersion $g$
and an embedding $h$, where $g$ and $h$ are also smooth maps.
   We merely let $g$ be the natural
projection $X \leftarrow X \times Y$ and $h$ be the embedding $X\times Y
\leftarrow Y$ of $Y$ onto the graph of $f$.   A similar factorization applies
in many other categories of mappings, 
but canonical relations require the following slightly
more complicated construction, since $g$ and $h$ above are not canonical
relations, even when $f$ is a symplectomorphism.  For any symplectic
manifold $Y$, 
$\epsilon_Y$ is the diagonal in $Y \times \ybar$, considered as 
a morphism to the  point $\unitob$ from $\ybar \times Y$.

\begin{prop}
\label{prop-factorization}
Any canonical relation $\xfy$ may be
factored into canonical relations as
$$ X = X \times \unitob
\stackrel{1_X\times \epsilon_Y}{\twoheadleftarrow} X \times\ybar \times Y
\stackrel{\gamma_f \times 1_Y}{\leftarrowtail} \unitob \times Y = Y. $$
This composition is strongly transversal, the relation $1_X\times
\epsilon_Y$ is a reduction, and $\gamma_f \times 1_Y$ is a coreduction. 
\end{prop}

\pf
In fact, the factorization works for any smooth relation (if we replace
$\ybar$ by $Y$ in the diagram).

We first look set-theoretically.   A pair $(x,y)$ lies in the
composition if and only if there is a triple $(x,y',y)\in X\times
\ybar\times Y$  for which $y'=y$ and $(x,y')\in f$, i.e. if and
only if $(x,y)\in f$.  Since this triple is determined by $x$ and $y$,
the composition is monic.   The composed relation, being
just $f$, is a closed submanifold, so we have proven monicity.

Next, we show that the composition is monic on the level of tangent
spaces.  In fact, we may simply apply the tangent functor to the
factorization above.   Identifying $T\unitob$ with $\unitob$, and
using the facts that $T1_X=1_{TX}$, $T\epsilon_Y = \epsilon_{TY}$,
$T\gamma_f = \gamma_{Tf}$, and $T1_Y=1_{TY},$
we obtain 
$$TX = TX \times\unitob
\stackrel{1_{TX}\times \epsilon_{TY}}{\twoheadleftarrow} TX \times
  T\ybar 
\times TY
\stackrel{\gamma_{Tf} \times 1_{TY}}{\leftarrowtail} \unitob \times TY
= TY. $$  But this is the same factorization as above, applied to $Tf$
rather than $f$, so it is monic, too.   

For canonical relations, monicity implies
transversality by symplectic duality.  In the general smooth case,
transversality may be verified directly; we omit the details.
 
Finally, $1_X\times \epsilon_Y$
is a reduction because $ 1_X$ and $\epsilon_Y$ are, and $\gamma_f
\times 1_Y$ is a coreduction because $\gamma_f$
and $1_Y$ are.   
\qed

\section{The Wehrheim-Woodward construction}

Wehrheim and Woodward \cite{we-wo:functoriality} construct a category
containing all the canonical relations between symplectic manifolds,
and in which composition coincides with set-theoretic composition in
the strongly transversal case.   

The construction of their category, which we will denote by 
$WW(\srel)$, begins with a category 
of ``paths''.   If we think of a category
as a directed graph with the objects as vertices and morphisms as
edges, these are paths in the usual sense, where we allow ``weakly
monotonic'' reparametrization.

\begin{dfn}
\label{dfn-paths}
The {\bf support} of an infinite composable sequence 
$$\boldlcf=(\ldots f_{-1},f_0,f_1,\ldots)$$
in any category 
is the set of integers $j$ for which $f_j$ is not an identity
morphism.  A {\bf canonical path in} $\srel$ is an infinite composable sequence
of {\em canonical} relations with finite support.  The target and source of $f_j$ for all
sufficiently large negative $j$ is thus a fixed object $X$ and, for
all sufficiently large positive $j$, a fixed object $Y$.  We call $X$
the {\bf target} and $Y$ the {\bf source} of the path $\boldlcf$.

Two canonical paths will be considered as equivalent if one may be obtained from
the other by inserting and removing finitely many identity morphisms.
This does not change the target or source.  The set of equivalence
classes is the {\bf path category} $P(\srel)$.   We will denote the
equivalence class of $(\ldots f_{-1},f_0,f_1,\ldots)$
by $ \langle \ldots f_{-1},f_0,f_1,\ldots \rangle$ and will also use the notation
$\langle f_r,\ldots,f_s \rangle $ when the support of the sequence $\boldlcf$ is contained
in the interval $[r,s]$.  

To compose $ \langle \boldlcf \rangle \in P(\srel)(X,Y)$ and $\langle \boldlcg
\rangle \in P(\srel)(Y,Z)$, 
choose representative sequences, remove
all but finitely many copies of $1_Y$ from the positive end of the
first sequence and the negative end of the second, and then concatenate the
truncated sequences. 

The identity morphism in $P(\srel)$ of any symplectic manifold $X$ 
is (represented by) the constant sequence 
with all entries equal to $1_X$.
\end{dfn}

\begin{rmk}
{\em  One could work as well with finite sequences, but the infinite version
is more convenient when it comes to defining rigid monoidal
structures.  
}
\end{rmk}

\begin{rmk}
\label{rmk-minimal}
{\em
Every morphism in $P(\srel)$ has a unique ``minimal'' representative for which
$f_i$ is an identity morphism for all $i\leq 0$ and for which there
are no identity morphisms in between nonidentity morphisms.  
}
\end{rmk}

\begin{rmk}
\label{rmk-shifting}
{\em 
A useful way to carry out the composition of two sequences is to shift the
first one (which does not change its equivalence class) 
so that its support (the set of $j$ such that $f_j$ is
not an identity morphism) is contained in the negative integers, and
to shift the second so that its support is contained in the positive
integers.  The composition is then represented by the sequence whose
value at $j$ is $f_j$ for $j\leq 0$ and $g_j$ for $j\geq 0$.

One may use a similar idea to verify associativity of composition;
given three sequences, shift them so that their supports are contained
in disjoint, successive intervals of integers.
}
\end{rmk}

We leave to the reader the proof of the following result.  
\begin{prop}
There is a unique functor $\srel \stackrel{c'}{\from} P(\srel)$
which is the identity on objects and which takes
each morphism $ \langle \ldots,f_{-1},f_0,f_1,\ldots \rangle $
to the composition $\cdots \circ f_{-1} \circ f_0 \circ f_1 \circ
\cdots $ in $\srel$.  (The ``infinite tails'' of identity morphisms
may be ignored here.)  
\end{prop}

We now define the Wehrheim-Woodward category $WW(\srel)$ by permitting  the 
actual composition of strongly transversal
pairs.   

\begin{dfn}
The Wehrheim-Woodward category
$WW(\srel)$ is the quotient category obtained from 
the category $P(\srel)$ of canonical paths by the
smallest equivalence relation for which two paths are equivalent
if a sequence representing one is obtained from 
a sequence representing the other by replacing successive entries
forming a strongly transversal pair $(p,q)$ by the single entry $pq$. 
The equivalence class in $WW(\srel)$ of $\langle \boldlcf \rangle \in
P(\srel\rangle$ will be denoted by $[\boldlcf]$. 
\end{dfn}

The composition functor $c'$ above descends to a 
functor $\srel \stackrel{c}{\from} WW(\srel)$, i.e.
$c([\ldots,f_{-1},f_{0},f_1,\ldots]) = \cdots \circ f_{-1}\circ
f_0\circ f_1 \circ \cdots.$

The canonical relations themselves
embed naturally in $WW(\srel)$ by the map $s$
(a cross section to the composition functor $c$) which maps each smooth morphism $f$ to the
equivalence class of sequences containing one entry equal to $f$ and
all the others equal to identity morphisms.
$WW(\srel)$ is then characterized by the universal property that 
any map which takes canonical relations to morphisms in some category $\boldc$,
which takes units to units, and which takes strongly transversal
compositions to compositions,
extends uniquely to a functor from $WW(\srel)$ to $\boldc$.

\section{Simplifying WW morphisms}
We prove here the main result of this note, that any morphism in
$WW(\srel)$ may be represented by a sequence of just two nontrivial
canonical relations.

\begin{thm}
\label{thm-two morphisms}
Let 
$(f_1,\ldots,f_r)$ be a composable sequence of canonical relations in $\srel$, with  
$f_i \in
\Hom(X_{i-1},X_i)$ 
for $i=1,\ldots r$. 
Then there is a symplectic manifold $Q$ with canonical relations 
$A \in \Hom(X_0,Q)$ and $B \in \Hom(Q,X_n)$ such that 
$A$ is a reduction, $B$ is a coreduction, and 
$[f_1,\ldots,f_r]=[A,B]$ in $WW(\srel)$.
\end{thm}

\pf
We illustrate the proof with diagrams for the case $r=4$, which is
completely representative of the general case.

First, we write $[f_{1},f_{2},f_{3},f_{4}]$ as a
  composition:

\[
\xymatrix{
X_0  && X_1  \ar[ll]_{f_{1}}  &&  X_2
\ar[ll]_{f_{2}} &&    X_3 \ar[ll]_{f_{3}}
&&  X_4 \ar[ll]_{f_{4}} 
}
\]

By Proposition \ref{prop-factorization}, we may factor each arrow as the
strongly transversal composition of a reduction and a coreduction.
The top row of the next
diagram is then equivalent to the zigzag line below it.

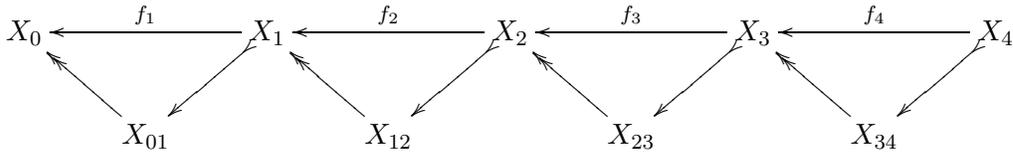
\[
\xymatrix{
X_0  && X_1  \ar@{>->}[dl] \ar[ll]_{f_{1}}  &&  X_2
\ar@{>->}[dl]\ar[ll]_{f_{2}} &&  \ar@{>->}[dl]  X_3 \ar[ll]_{f_{3}}
&&  X_4 \ar@{>->}[dl]\ar[ll]_{f_{4}} 
\\
  & X_{01}  \ar@{->>}[ul]&& X_{12}  \ar@{->>}[ul]&& 
            X_{23} \ar@{->>}[ul] 
      && X_{34} \ar@{->>}[ul]\\
 }
\]

Next, we compose pairs of diagonal arrows to produce the bottom row below.

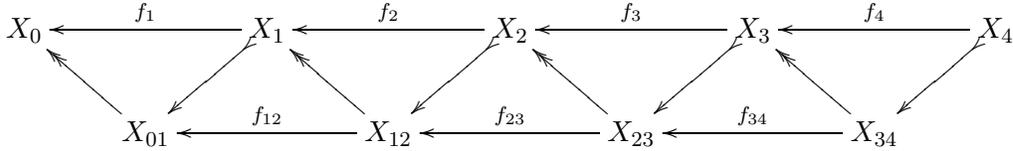
\[
\xymatrix{
X_0  && X_1  \ar@{>->}[dl]\ar[ll]_{f_{1}}  &&  X_2
\ar@{>->}[dl]\ar[ll]_{f_{2}} &&  \ar@{>->}[dl]  X_3 \ar[ll]_{f_{3}}
&&  X_4 \ar@{>->}[dl]\ar[ll]_{f_{4}} 
\\
  & X_{01}  \ar@{->>}[ul]&& X_{12}  \ar@{->>}[ul]\ar[ll] _{f_{12}}&& 
            X_{23} \ar@{->>}[ul]\ar[ll] _{f_{23}} 
      && X_{34} \ar@{->>}[ul]\ar[ll]  _{f_{34}} \\
 }
\]

Each of these compositions is strongly transversal, even ``doubly so'', thanks
to the decorations on the arrows identifying them as reductions and
coreductions.  It follows that the original composition on the top row
is equivalent to the composition of the bottom row with the outer
diagonal edges.   

We repeat the process to obtain another row.

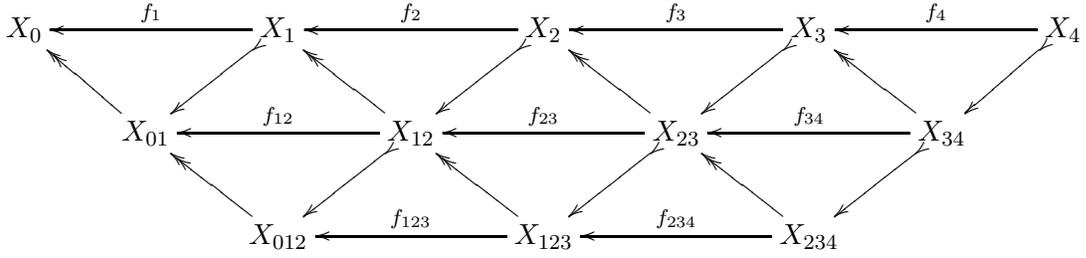
\[
\xymatrix{
X_0  && X_1  \ar@{>->}[dl]\ar[ll]_{f_{1}}  &&  X_2
\ar@{>->}[dl]\ar[ll]_{f_{2}} &&  \ar@{>->}[dl]  X_3 \ar[ll]_{f_{3}}
&&  X_4 \ar@{>->}[dl]\ar[ll]_{f_{4}} 
\\
  & X_{01}  \ar@{->>}[ul]&& X_{12} \ar@{>->}[dl] \ar@{->>}[ul]\ar[ll] _{f_{12}}&& 
            X_{23}  \ar@{>->}[dl]\ar@{->>}[ul]\ar[ll] _{f_{23}} 
      && X_{34} \ar@{->>}[ul]\ar@{>->}[dl]\ar[ll]  _{f_{34}}\\
   && X_{012}  \ar@{->>}[ul]&& X_{123} \ar@{->>}[ul]\ar[ll] _{f_{123}} &&
   X_{234} \ar@{->>}[ul]\ar[ll] _{f_{234}}
}
\]

Repeating two more times, we arrive at a triangle, in which the top
row is equivalent to the composition of the arrows on the other two
sides.  

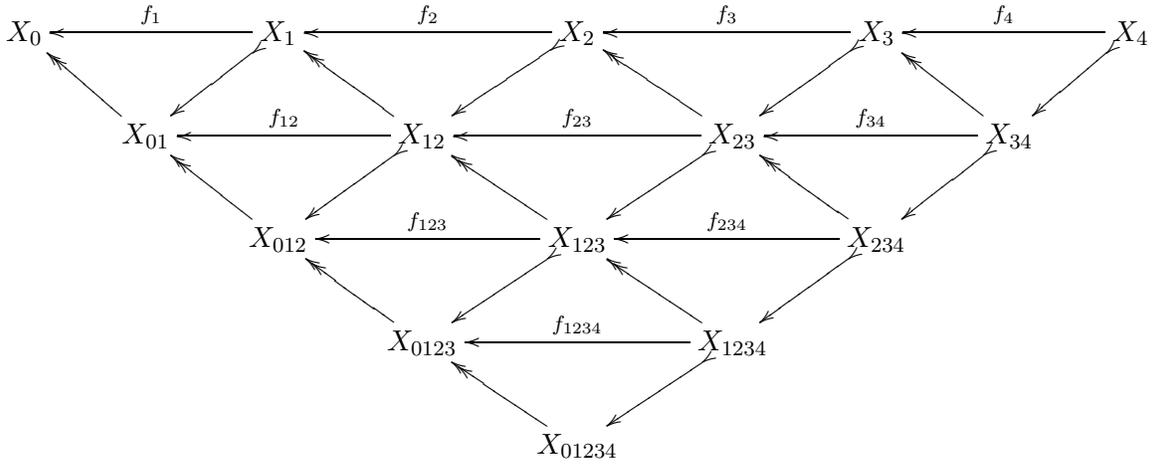
\[
\xymatrix{
X_0  && X_1  \ar@{>->}[dl]\ar[ll]_{f_{1}}  &&  X_2
\ar@{>->}[dl]\ar[ll]_{f_{2}} &&  \ar@{>->}[dl]  X_3 \ar[ll]_{f_{3}}
&&  X_4 \ar@{>->}[dl]\ar[ll]_{f_{4}} 
\\
  & X_{01}  \ar@{->>}[ul]&& X_{12} \ar@{>->}[dl] \ar@{->>}[ul]\ar[ll] _{f_{12}}&& 
            X_{23}  \ar@{>->}[dl]\ar@{->>}[ul]\ar[ll] _{f_{23}} 
      && X_{34} \ar@{->>}[ul]\ar@{>->}[dl]\ar[ll]  _{f_{34}}\\
   && X_{012}  \ar@{->>}[ul]&& X_{123} \ar@{>->}[dl] \ar@{->>}[ul]\ar[ll] _{f_{123}} &&
   X_{234} \ar@{->>}[ul]\ar @{>->}[dl]\ar[ll] _{f_{234}}
\\
     &&&  X_{0123} \ar@{->>}[ul]&& X_{1234}\ar @{>->}[dl] \ar@{->>}[ul]\ar[ll] _{f_{1234}}\\
      &&&&     
X_{01234}  \ar @{->>}[ul]
}
\]

Finally, we observe that all the arrows going up the left-hand side
are reductions, so we may compose them all to produce a single
reduction $A$.  Similarly, the coreductions going down on the right
yield a coreduction $B$.

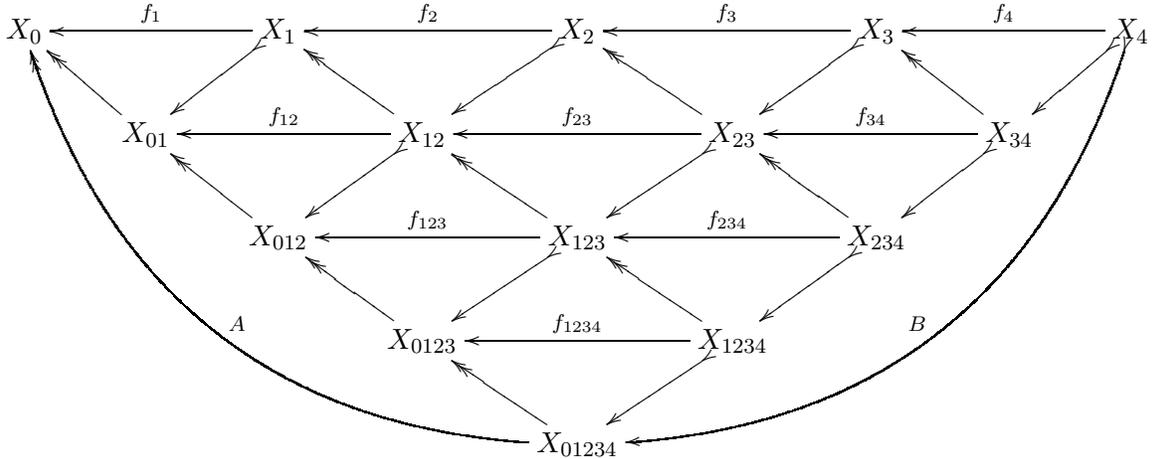
\[
\xymatrix{
X_0  && X_1  \ar@{>->}[dl]\ar[ll]_{f_{1}}  &&  X_2
\ar@{>->}[dl]\ar[ll]_{f_{2}} &&  \ar@{>->}[dl]  X_3 \ar[ll]_{f_{3}}
&&  X_4 \ar@{>->}[dl]\ar[ll]_{f_{4}} 
\ar@/^4pc/@{>->}[ddddllll]_B
\\
  & X_{01}  \ar@{->>}[ul]&& X_{12} \ar@{>->}[dl] \ar@{->>}[ul]\ar[ll] _{f_{12}}&& 
            X_{23}  \ar@{>->}[dl]\ar@{->>}[ul]\ar[ll] _{f_{23}} 
      && X_{34} \ar@{->>}[ul]\ar@{>->}[dl]\ar[ll]  _{f_{34}}\\
   && X_{012}  \ar@{->>}[ul]&& X_{123} \ar@{>->}[dl] \ar@{->>}[ul]\ar[ll] _{f_{123}} &&
   X_{234} \ar@{->>}[ul]\ar @{>->}[dl]\ar[ll] _{f_{234}}
\\
     &&&  X_{0123} \ar@{->>}[ul]&& X_{1234}\ar @{>->}[dl] \ar@{->>}[ul]\ar[ll] _{f_{1234}}\\
      &&&&  \ar @/^4pc/@{->>}[uuuullll]_A 
X_{01234}  \ar @{->>}[ul]
}
\]

We may now erase everything in the middle of the diagram to obtain the desired
factorization.

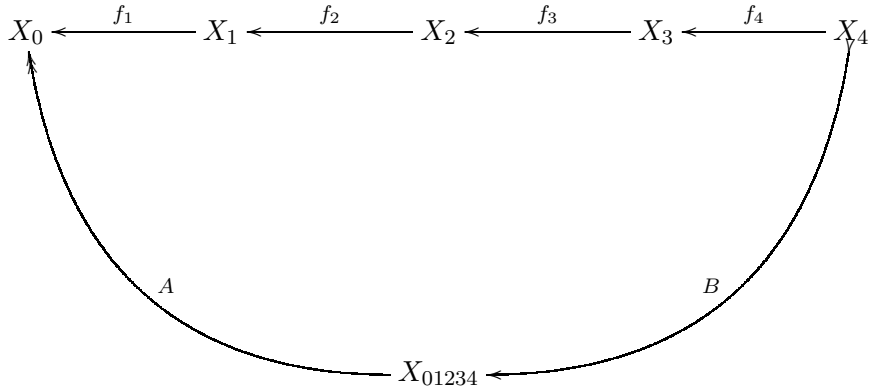
\[
\xymatrix{
X_0  && X_1 \ar[ll]_{f_{1}}  &&  X_2
\ar[ll]_{f_{2}} &&    X_3 \ar[ll]_{f_{3}}
&&  X_4 \ar[ll]_{f_{4}} 
\ar@/^4pc/@{>->}[ddddllll]_B
\\
 \\
\\
   \\
      &&&&  \ar @/^4pc/@{->>}[uuuullll]_A 
X_{01234} 
}
\]
\qed

\begin{rmk}
{\em Sylvain Cappell has pointed out the similarity of this result to ideas
of J.H.C. Whitehead on simple homotopy theory, where maps are factored
in to collapses and expansions \cite{co:course}.   And Thomas Kragh has noted a resemblance
to the theory of Waldhausen categories; the diagram on page 207 of
\cite{ro:lecture} looks very much like the ones in the proof above.
}\end{rmk}

\end{document}